# Mô hình toán mô tả quỹ đạo chuyển động của trục trong ổ trục của bơm bánh răng ăn khớp trong

**Nguyễn Mạnh Hùng**
*Bộ môn Đại số – Xác suất thống kê, Khoa Khoa học Cơ bản, ĐH GTVT*
*Email: hung.manh.nguyen@utc.edu.vn*

**Phạm Trọng Hòa**
*Bộ môn Máy xây dựng, Khoa Cơ Khí, ĐH GTVT*
*Email: hoagtvt100@gmail.com*

***Tóm tắt***: *Bài báo này trình bày về mô hình toán xác định chuyển động của trục trong ổ trục của bơm bánh răng ăn khớp trong. Khi bơm làm việc, trục bánh răng dao động xung quanh vị trí cân bằng dưới tác động của lớp màng dầu bôi trơn giữa trục và thân bơm. Để xác định quỹ đạo chuyển động của trục, lớp màng dầu bôi trơn được mô hình hóa thành một hệ lò xo và giảm xóc. Các hệ số độ cứng và hệ số giảm xóc được tính dựa vào phân bố áp suất của lớp màng dầu lên trục. Theo lý thuyết bôi trơn, phân bố áp suất của lớp màng dầu có thể tìm được bằng cách giải phương trình Reynold. Cuối cùng, chuyển động của trục bánh răng được xác định bởi dao động tắt dần của hệ lò xo và giảm xóc. Bài báo cũng thực hiện tính toán mô phỏng tìm quỹ đạo chuyển động của trục với các dữ kiện đầu vào cho trước. Các kết quả tính toán có thể được sử dụng để xác định điều kiện làm việc tối ưu của bơm, tránh các trường hợp tiếp xúc giữa trục và thân bơm gây ăn mòn vật liệu.*

## I. Giới thiệu

Màng dầu bôi trơn thường được sử dụng trong các thiết bị máy móc để ngăn cách và làm trơn các bề mặt tiếp xúc. Đối với các máy chuyển động quay, lớp màng dầu bôi trơn có vai trò quan trọng trong việc truyền tải các lực tác động thông qua áp suất sinh ra trong lớp màng dầu. Do đó tính chất hoạt động của máy móc có thể phân tích được dựa vào các quy luật thủy động lực học. Các vấn đề cơ bản của lý thuyết bôi trơn có thể được tìm hiểu chi tiết trong [1].

Nhiều nghiên cứu cho đến nay đã được thực hiện để xác định các tính chất thủy động học của các loại ổ đỡ khác nhau. Chẳng hạn như trong [2], Z. L. Qiu xây dựng mô hình để tính toán các đặc trưng tĩnh và động của một loại ổ đỡ, đồng thời phân tích tính chất ổn định của nó. Các kết quả thí nghiệm cho thấy sự phù hợp của các kết quả mô phỏng. P. T. Hòa [3-5] nghiên cứu vấn đề chuyển động của ổ trục của bơm bánh răng ăn khớp trong bằng phương pháp mobility, được đề xuất lần đầu năm 1965 bởi J. F. Booker [6]. Tính ổn định của hệ bánh răng loại này cũng được xem xét trong [7, 8]. Nói chung, trong phân tích các đặc trưng chuyển động và ổn định của ổ trục, lực tác động của màng dầu được tuyến tính hóa và biểu diễn như là một hàm của chuyển vị và vận tốc của trục quay (đề xuất bởi J. W. Lund [9]). Sự tuyến tính hóa này dẫn đến việc xem xét lớp màng dầu như một hệ lò xo và giảm xóc. Kết quả tính toán các hệ số độ cứng và hệ số giảm xóc cho phép phân tích các đặc trưng động học của ổ trục như trong trường hợp bơm bánh răng ăn khớp trong (xem P. T. Hòa [10]).

Cho đến nay, chưa có một mô hình toán nào được đề xuất giúp xác định quỹ đạo chuyển động của ổ trục của bơm bánh răng ăn khớp trong. Một mô hình tương tự như thế đã được xây dựng bởi R. Castilla [11] để phân tích chuyển động của trục quay của một loại bơm bánh răng; và kết quả tính toán cũng tương thích với số liệu thực nghiệm. Trong báo cáo này, chúng tôi nghiên cứu một mô hình tính toán quỹ đạo chuyển động cho ổ trục của bơm bánh răng ăn khớp trong. Mô hình này, sau khi được kiểm nghiệm với dữ liệu thực tế, có thể được sử dụng để xác định điều kiện làm việc tối ưu cho bơm.

## II. Mô hình toán mô tả chuyển động

**Hình 1** mô tả cấu tạo của bơm bánh răng ăn khớp trong, ở đó lớp màng dầu ngăn cách vành bánh răng và thân bơm tiếp xúc. Đồng thời lớp màng dầu cũng có vai trò trong việc truyền lực tác động trong quá trình hoạt động thông qua áp suất trong lớp màng dầu.

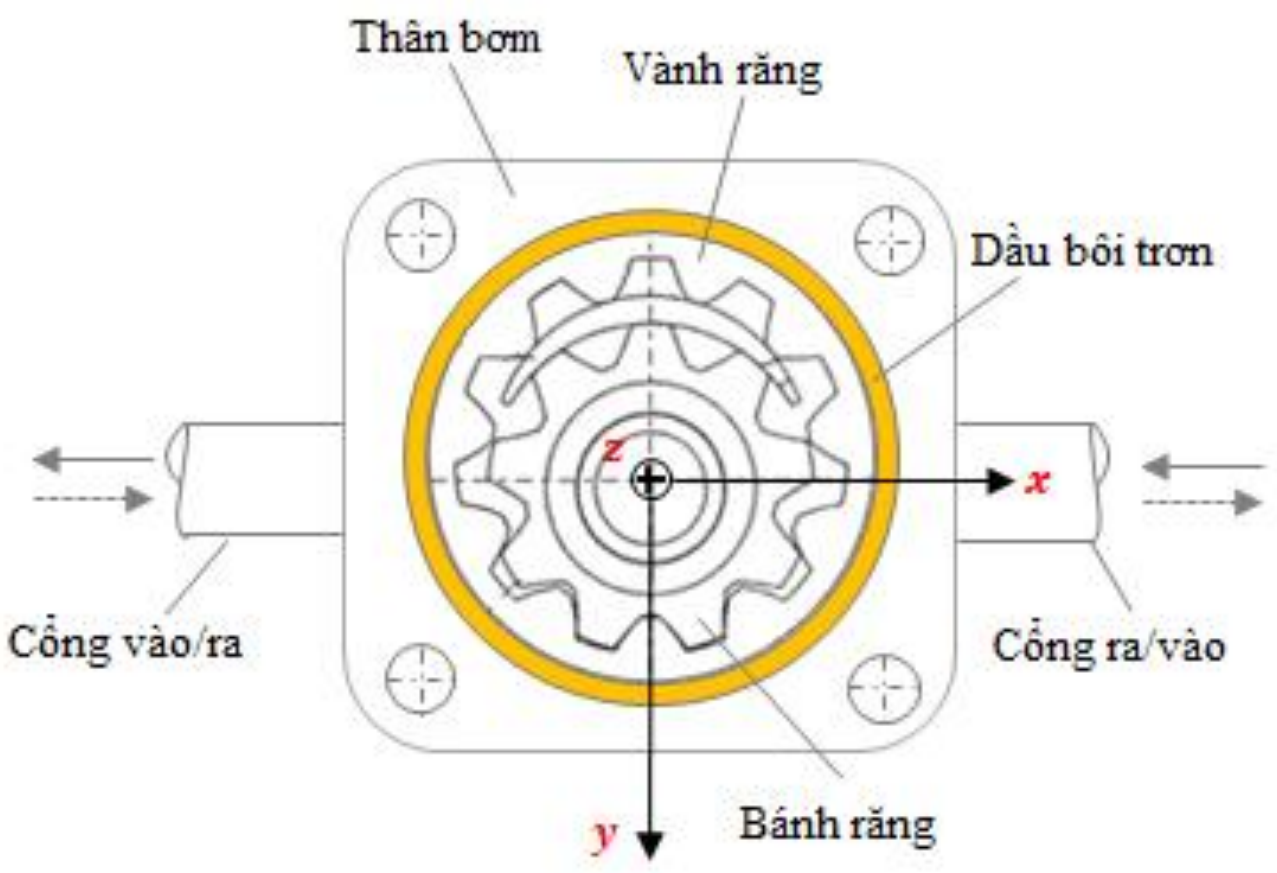

**Hình 1**: Mặt cắt của bơm bánh răng ăn khớp trong

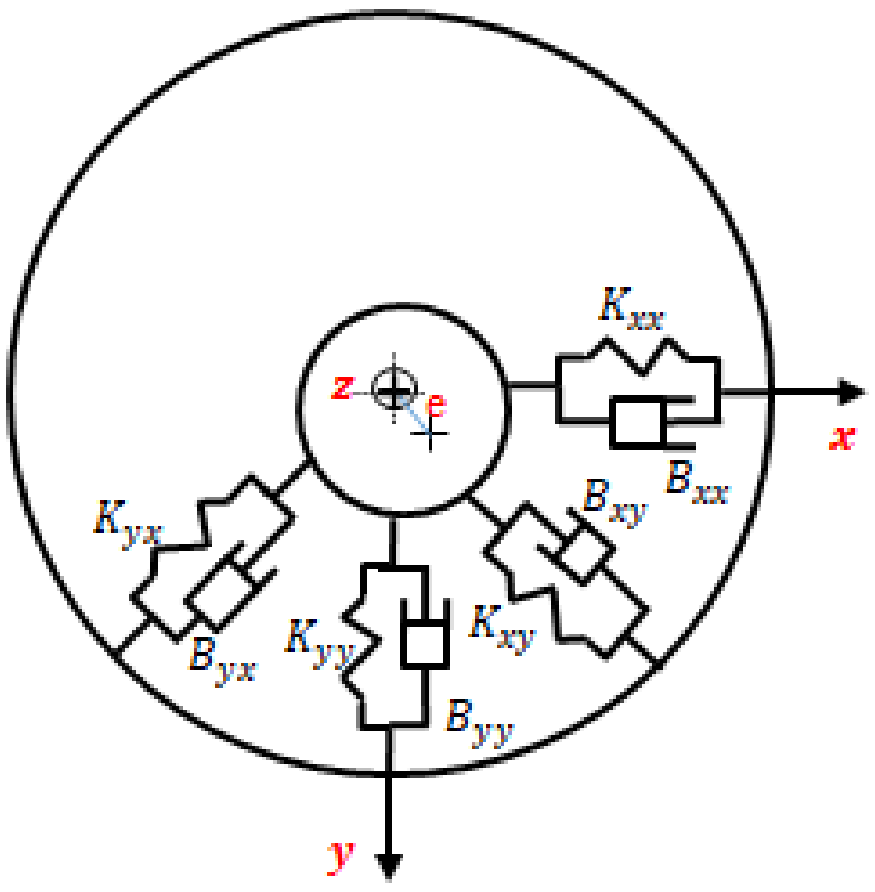

**Hình 2**: Mô hình lò xo và giảm xóc

Trong phân tích chuyển động của hệ bánh răng, một phương pháp đơn giản có thể thực hiện được là mô hình hóa lớp màng dầu như một hệ lò xo và giảm xóc như ở **Hình 2**. Hệ thống này được mô tả bởi 4 hệ số độ cứng $k_{xx}$, $k_{xy}$, $k_{yx}$, $k_{yy}$ và 4 hệ số giảm xóc $b_{xx}$, $b_{xy}$, $b_{yx}$, $b_{yy}$. Tám hệ số này xác định các đặc trưng động học của hệ bánh răng, và do đó được gọi là các hệ số động học. Các hệ số này là hàm của tỉ số độ lệch tâm ε (e/c), trong đó c được gọi là khe hở xuyên tâm giữa vành răng và thân bơm, và e là độ lệch tâm của vành răng. Phương pháp mô hình hóa này khá phù hợp khi trục bánh răng có chuyển động tương đối nhỏ xung quanh vị trí cân bằng.

### 1. Phương trình Reynold

Xét hệ trục bánh răng và thân bơm được mô tả trong **Hình 3**. Giả thiết rằng khoảng trống giữa trục bánh răng và thân bơm được lấp đầy dầu bôi trơn. Hệ trục tọa độ Oxy có gốc tọa độ là tâm của thân bơm.

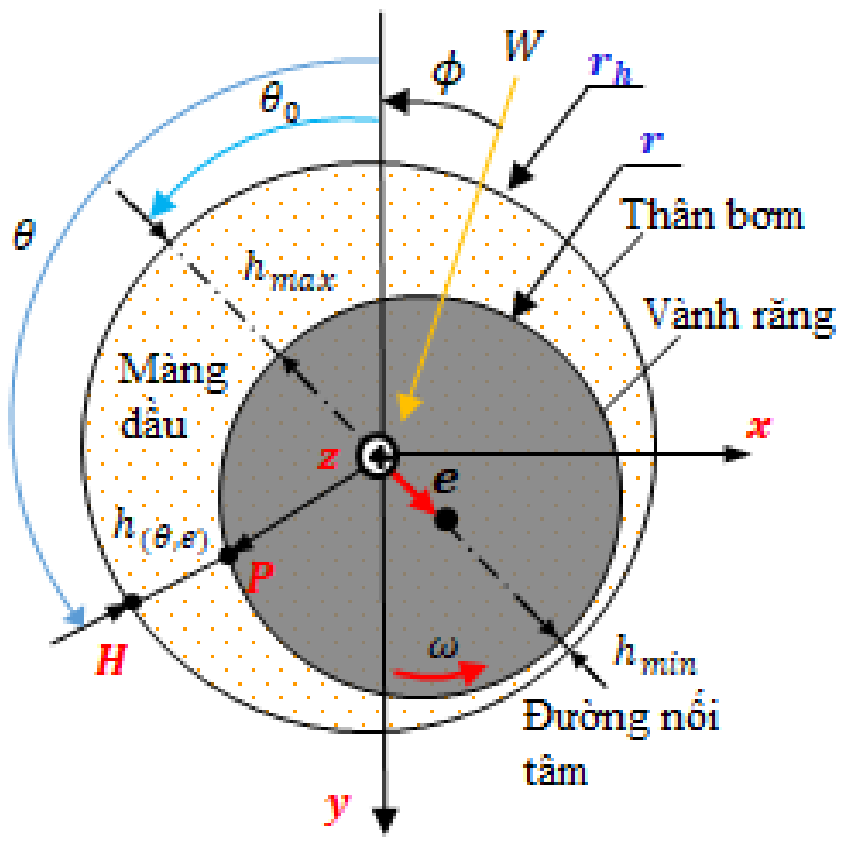

**Hình 3**: Mô tả hình học lớp màng dầu và vành bánh răng

Phân bố áp suất $p$ trong lớp màng dầu được xác định bởi phương trình Reynold cho dòng chảy tầng dưới đây (xem [9]):

$$\frac{1}{r^2}\frac{\partial}{\partial\theta}\left(\frac{\rho h^3}{12\mu}\frac{\partial p}{\partial\theta}\right)+\frac{\partial}{\partial z}\left(\frac{\rho h^3}{12\mu}\frac{\partial p}{\partial z}\right)=\frac{\omega}{2}\frac{\partial(\rho h)}{\partial\theta}+\rho\left(V_x\sin\theta+V_y\cos\theta\right)+\dot{\rho}h \qquad (1)$$

Theo [1], độ dày lớp màng dầu xấp xỉ bằng:

$$h=c+e\cos(\theta-\theta_0)=c+x.\sin\theta+y.\cos\theta$$

Giả thiết rằng dòng chất lưu có độ nhớt đồng nhất và không chịu nén thì phương trình Reynold có thể được viết lại dưới dạng sau:

$$\frac{1}{r^2}\frac{\partial}{\partial\theta}\left(\frac{h^3}{12\mu}\frac{\partial p}{\partial\theta}\right)+\frac{\partial}{\partial z}\left(\frac{h^3}{12\mu}\frac{\partial p}{\partial z}\right)=\frac{\omega}{2}\frac{\partial h}{\partial\theta}+\dot{x}\sin\theta+\dot{y}\cos\theta \qquad (2)$$

Đạo hàm cả hai vế phương trình (2) theo $x, y, \dot{x}$ và $\dot{y}$, và kết hợp với phương trình (2) ở trạng thái tĩnh (cân bằng), tức là $\dot{x}=\dot{y}=0$, ta thu được hệ 5 phương trình:

$$\frac{1}{r^2}\frac{\partial}{\partial\theta}\left(\frac{h^3}{12\mu}\frac{\partial}{\partial\theta}\right)+\frac{\partial}{\partial z}\left(\frac{h^3}{12\mu}\frac{\partial}{\partial z}\right)\begin{pmatrix}p_0\\p_x\\p_y\\p_{\dot{x}}\\p_{\dot{y}}\end{pmatrix}=\begin{pmatrix}\frac{\omega}{2}\frac{\partial h}{\partial\theta}\\ \frac{\omega}{2}\left(\cos\theta-\frac{3\sin\theta}{h}\frac{\partial h}{\partial\theta}\right)-\frac{h^3}{4\mu r^2}\frac{\partial}{\partial\theta}\left(\frac{\sin\theta}{h}\right)\frac{\partial p_0}{\partial\theta}\\ -\frac{\omega}{2}\left(\sin\theta+\frac{3\cos\theta}{h}\frac{\partial h}{\partial\theta}\right)-\frac{h^3}{4\mu r^2}\frac{\partial}{\partial\theta}\left(\frac{\cos\theta}{h}\right)\frac{\partial p_0}{\partial\theta}\\ \sin\theta\\ \cos\theta\end{pmatrix} \qquad (3)$$

trong đó $p_0$ là áp suất tĩnh của lớp màng dầu và

$$p_x=\left(\frac{\partial p}{\partial x}\right)_0,\qquad p_y=\left(\frac{\partial p}{\partial y}\right)_0,\qquad p_{\dot{x}}=\left(\frac{\partial p}{\partial\dot{x}}\right)_0,\qquad p_{\dot{y}}=\left(\frac{\partial p}{\partial\dot{y}}\right)_0$$

Để giải hệ phương trình (3), đầu tiên ta tìm áp suất tĩnh $p_0$ bằng lược đồ sai phân với lưới chia như trong **Hình 4b**.

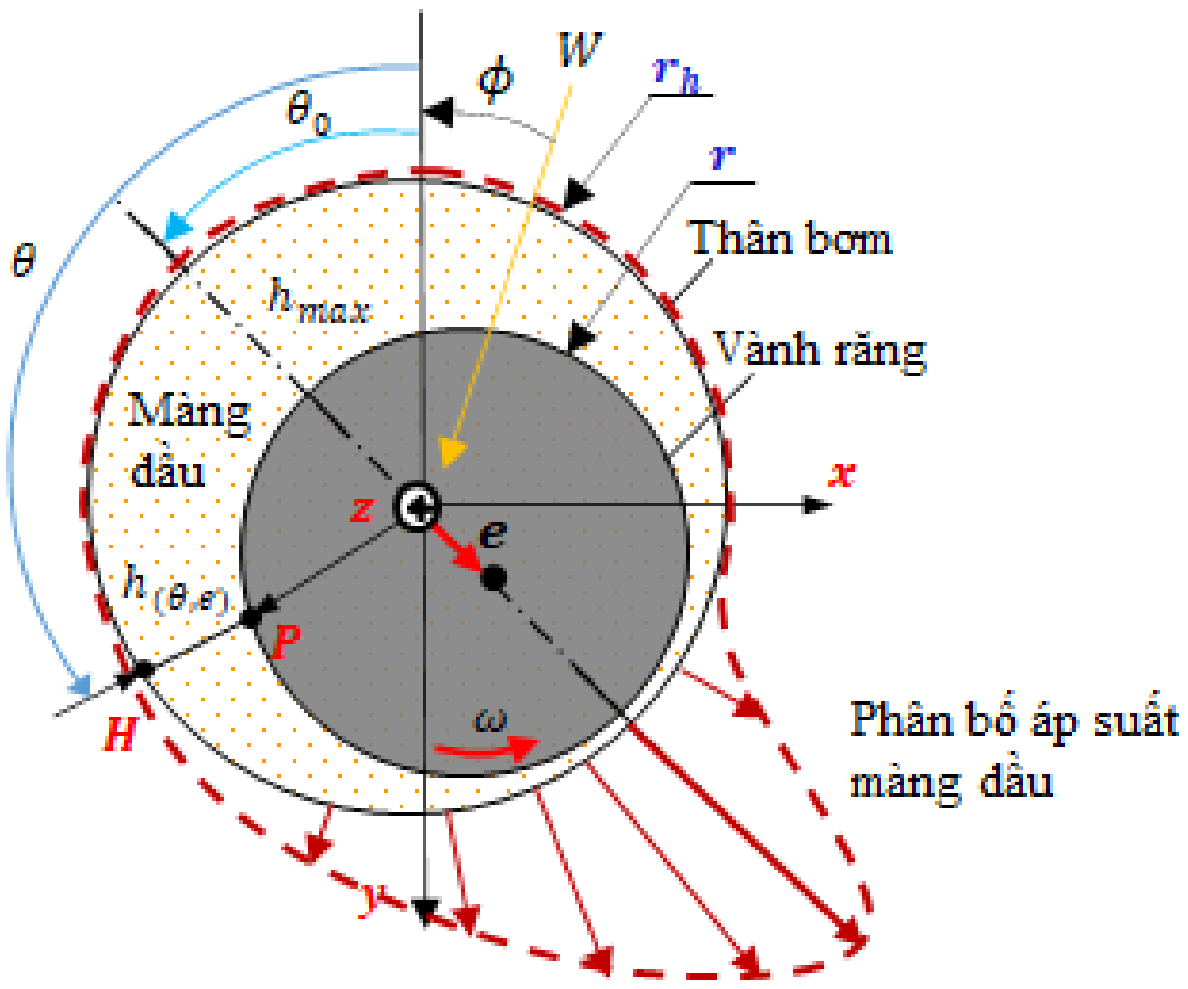

**Hình 4a**: Phân bố áp suất màng dầu

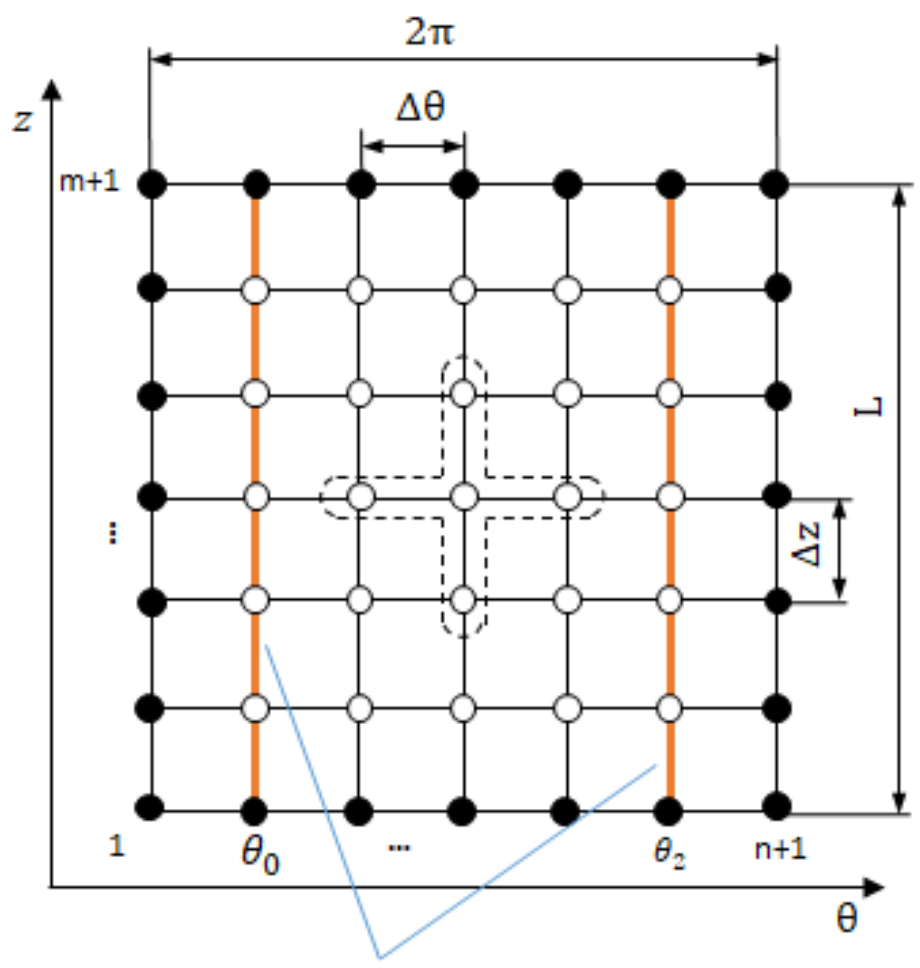

**Hình 4b**: Lưới sai phân

Các thí nghiệm nhằm xác định phân bố lớp màng dầu cho thấy áp suất dương chỉ phân bố trên một phần bề mặt của vành răng giống như mô tả trên hình vẽ. Do đó, điều kiện biên Reynold thường được chọn thỏa mãn:

- $p(\theta, 0) = p(\theta, L) = 0$
- $p(\theta_0, z) = p(\theta_2, z) = 0$

Biên Reynolds $\theta_0$ và $\theta_2$ không biết trước và sẽ được xác định trong quá trình tính toán. Để giải hệ phương trình Reynold, $\theta_0$ và $\theta_2$ được gán cho các giá trị xuất phát, sau đó giải ra $p_0$ rồi điều chỉnh $\theta_0$ và $\theta_2$ bằng cách gán các giá trị áp suất âm bằng 0 (xem [2]). Lặp đi lặp lại quá trình tính toán đó cho đến khi hội tụ thì dừng. Kết quả nhận được là biên $\theta_0,\ \theta_2$ và phân bố áp suất tĩnh $p_0$. Sử dụng các kết quả này và giải bốn phương trình còn lại bằng lược đồ sai phân để tìm $p_x$ , $p_y$ , $p_{\dot{x}}$ và $p_{\dot{y}}$.

## 2. Lực màng dầu tác động lên hệ bánh răng

Lực tác động của màng dầu được gây nên bởi áp suất màng dầu được xác định bởi công thức:

$$\begin{bmatrix} F_x \\ F_y \end{bmatrix} = -\int_{\theta_0}^{\theta_2}\int_0^L p \begin{bmatrix} \sin\theta \\ \cos\theta \end{bmatrix} (r d\theta)\, dz$$

Theo phương pháp của Lund [9], lực tác động này là hàm của vị trí và vận tốc của trục bánh răng và được tuyến tính hóa dưới dạng sau:

$$\begin{bmatrix} F_x \\ F_y \end{bmatrix} = \begin{bmatrix} (F_x)_0 \\ \left(F_y\right)_0 \end{bmatrix} + \begin{bmatrix} k_{xx} & k_{xy} \\ k_{yx} & k_{yy} \end{bmatrix}\begin{bmatrix} \Delta x \\ \Delta y \end{bmatrix} + \begin{bmatrix} b_{xx} & b_{xy} \\ b_{yx} & b_{yy} \end{bmatrix}\begin{bmatrix} \Delta \dot{x} \\ \Delta \dot{y} \end{bmatrix} \tag{4}$$

trong đó $(F_x)_0$ và $\left(F_y\right)_0$ là các lực gây ra bởi áp suất tĩnh của màng dầu ở trạng thái cân bằng:

$$\begin{bmatrix} (F_x)_0 \\ \left(F_y\right)_0 \end{bmatrix} = -r\int_{\theta_0}^{\theta_2}\int_0^L p_0 \begin{bmatrix} \sin\theta \\ \cos\theta \end{bmatrix} d\theta\, dz \tag{5}$$

Các hệ số độ cứng và hệ số giảm xóc được xác định bởi công thức:

$$\begin{bmatrix} k_{xx} & k_{xy} \\ k_{yx} & k_{yy} \end{bmatrix} = \begin{bmatrix} \dfrac{\partial F_x}{\partial x} & \dfrac{\partial F_x}{\partial y} \\ \dfrac{\partial F_y}{\partial x} & \dfrac{\partial F_y}{\partial y} \end{bmatrix} = -r\int_{\theta_0}^{\theta_2}\int_0^L \begin{bmatrix} p_x \sin\theta & p_y \sin\theta \\ p_x \cos\theta & p_y \cos\theta \end{bmatrix} d\theta\, dz \tag{6}$$

$$\begin{bmatrix} b_{xx} & b_{xy} \\ b_{yx} & b_{yy} \end{bmatrix} = \begin{bmatrix} \dfrac{\partial F_x}{\partial \dot{x}} & \dfrac{\partial F_x}{\partial \dot{y}} \\ \dfrac{\partial F_y}{\partial \dot{x}} & \dfrac{\partial F_y}{\partial \dot{y}} \end{bmatrix} = -r\int_{\theta_0}^{\theta_2}\int_0^L \begin{bmatrix} p_{\dot{x}} \sin\theta & p_{\dot{y}} \sin\theta \\ p_{\dot{x}} \cos\theta & p_{\dot{y}} \cos\theta \end{bmatrix} d\theta\, dz \tag{7}$$

Các số hạng trên có thể tính được bằng cách sử dụng phương pháp Simpson tính gần đúng tích phân.

## 3. Phương trình chuyển động của vành răng

Ký hiệu $m_a$ là khối lượng của hệ trục bánh răng và $W$ là ngoại lực tác động vào hệ theo phương tạo với trục thẳng đứng một góc $\phi$. Phương trình chuyển động của trục bánh răng được viết như sau:

$$m_a\begin{bmatrix}\Delta\ddot{x}\\ \Delta\ddot{y}\end{bmatrix}=\begin{bmatrix}W\sin\phi\\ W\cos\phi\end{bmatrix}-\begin{bmatrix}(F_x)_0\\ (F_y)_0\end{bmatrix}-\begin{bmatrix}k_{xx} & k_{xy}\\ k_{yx} & k_{yy}\end{bmatrix}\begin{bmatrix}\Delta x\\ \Delta y\end{bmatrix}-\begin{bmatrix}b_{xx} & b_{xy}\\ b_{yx} & b_{yy}\end{bmatrix}\begin{bmatrix}\Delta\dot{x}\\ \Delta\dot{y}\end{bmatrix} \tag{8}$$

Lúc bắt đầu tính toán quỹ đạo chuyển động, trục bánh răng được giả thiết nằm ở một vị trí xác định, tức cho trước $e$ và $\theta_0$. Giải các phương trình Reynold (1) để tìm áp suất tĩnh và các đạo hàm của nó. Các hệ số động học được tính toán theo các công thức (3) – (5). Sau đó giải phương trình (6) bằng phương pháp Runge – Kutta bậc 4 và cập nhật lại trạng thái $e$ và $\theta_0$ của trục bánh răng. Lặp lại quá trình tính toán trên cho thời điểm tiếp theo, ta sẽ thu được quỹ đạo chuyển động của hệ bánh răng.

## III. Các kết quả tính toán số

### *1. Sơ đồ tính toán*

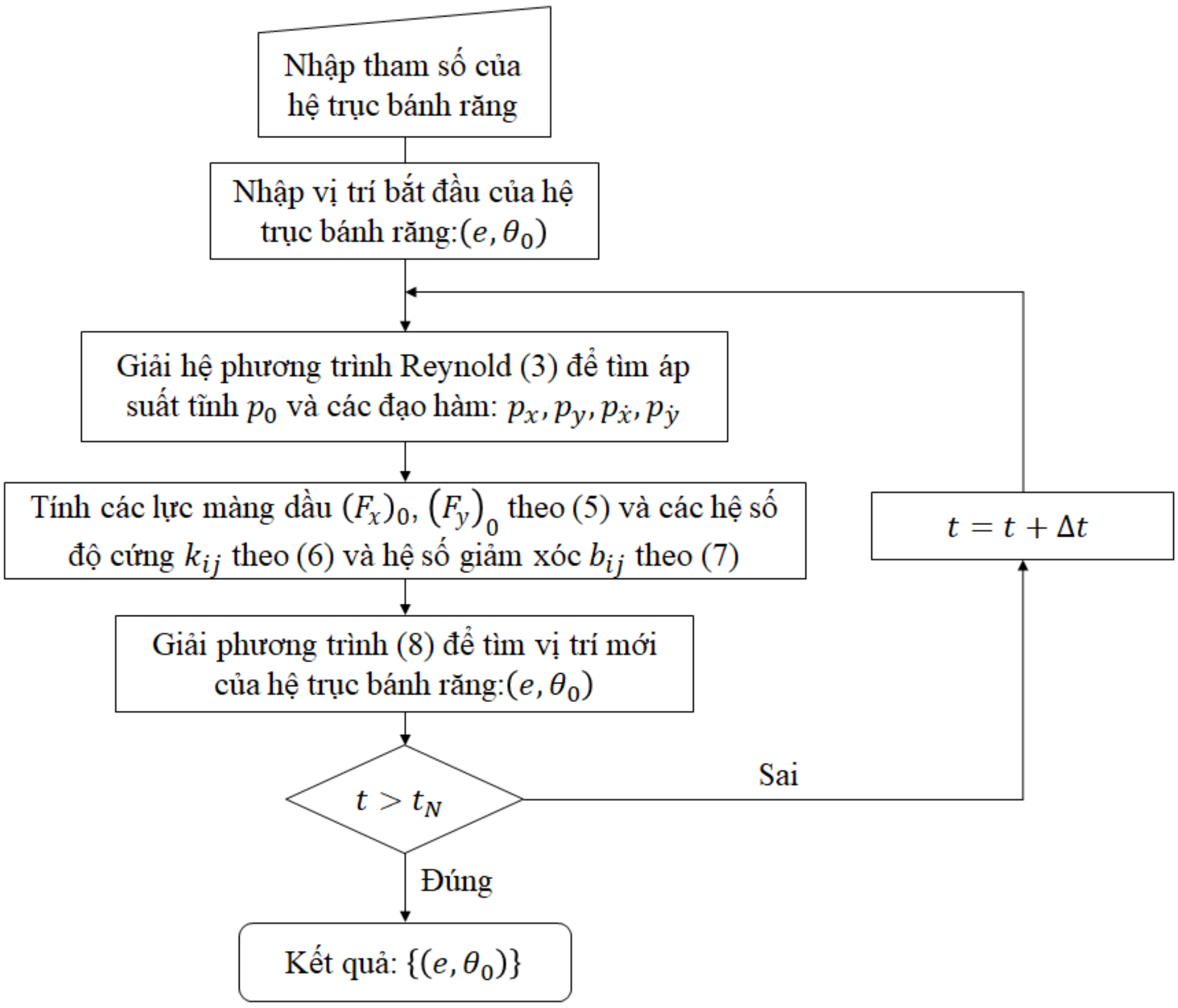

**Hình 5**: Sơ đồ tính toán quỹ đạo chuyển động

### *2. Một số kết quả mô phỏng*

Sau đây chúng ta tiến hành mô phỏng chuyển động của trục bánh răng của bơm bánh răng ăn khớp trong với các tham số đầu vào được cho dưới đây:

**Bảng 1**: Các tham số đầu vào của quá trình mô phỏng

| **Tham số** | **Ký hiệu** | **Giá trị** | **Đơn vị** |
|---|---|---|---|
| Bán kính vành răng | r | 0.055 | m |
| Chiều dài của trục | L | 0.034 | m |
| Độ nhớt | $\mu$ (HLP46) | 0.041 | Pa.s |
| Vận tốc quay | $\omega$ | 3000 | Vòng/phút |
| Khe hở xuyên tâm | c | $5.10^{-5}$ | m |

Mô phỏng được thực hiện với các tải trọng tác dụng theo phương thẳng đứng với giá trị khác nhau, từ 1000 N đến 4000 N, và với các vận tốc quay của trục bánh răng thay đổi từ 1000 vòng/phút cho đến 3000 vòng/phút. Các kết quả tính toán được thể hiện như dưới đây.

*a) Quỹ đạo của trục bánh răng dưới tác động của tải trọng cố định:*

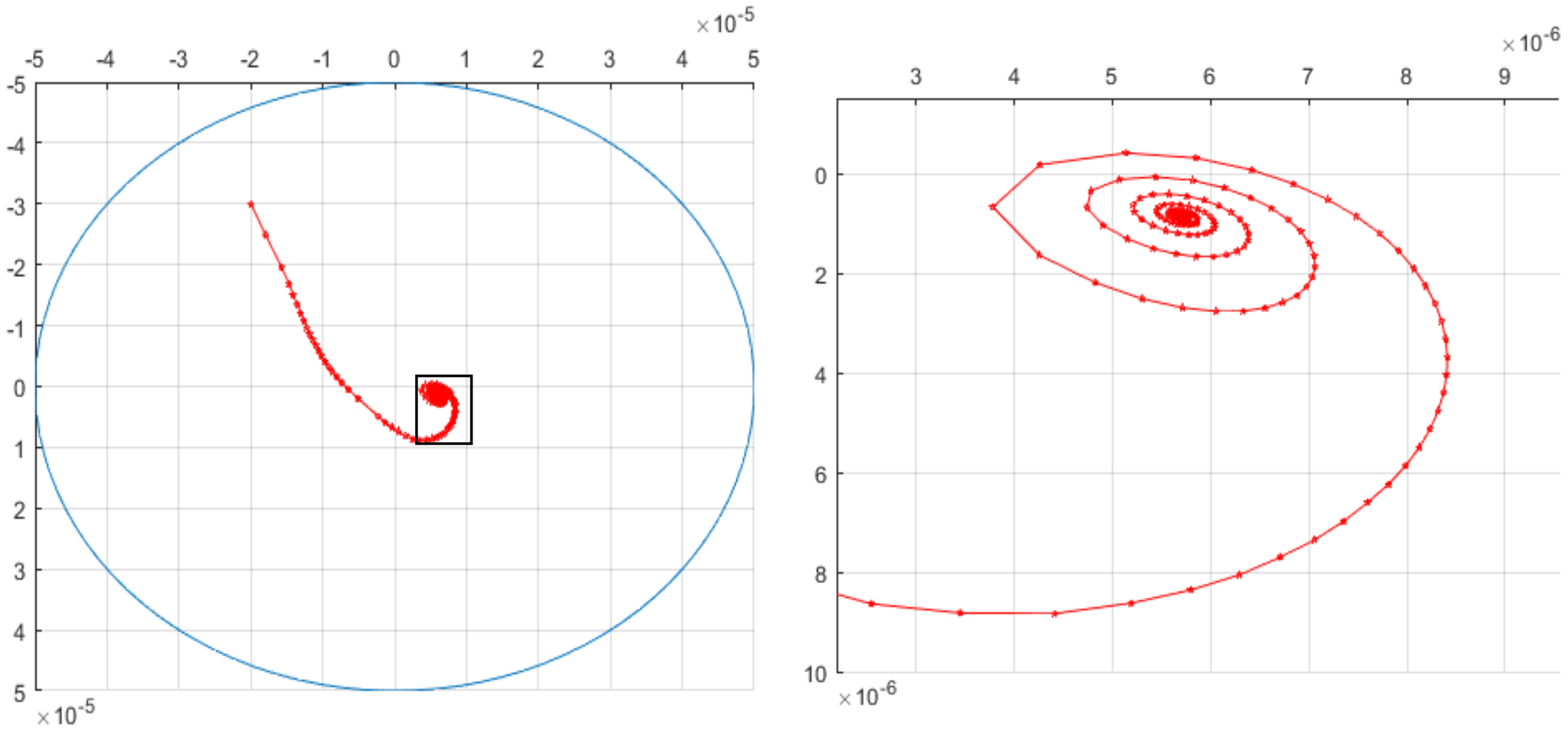

**Hình 5**: Quỹ đạo chuyển động (trái) và dao động (phải) với $W = 1000$ (N)

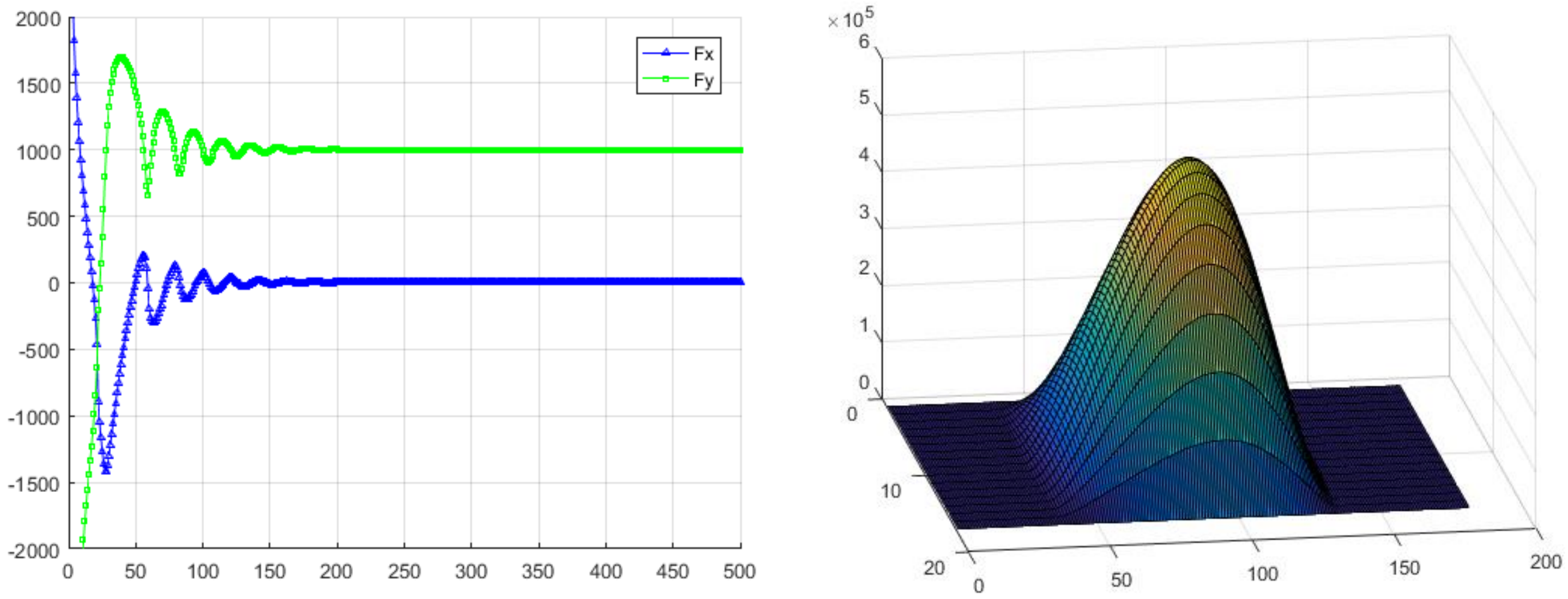

**Hình 6**: Lực màng dầu (trái) và phân bố áp suất tĩnh (phải) với $W = 1000$ (N)

**Hình 5** (trái) mô tả quỹ đạo của trục bơm bánh răng ăn khớp trong khi tâm trục xuất phát tại vị trí có tọa độ $(-2.10^{-5}\text{m};-3.10^{-5}\text{m})$ và chịu tải trọng $W = 1000$ (N) theo phương thẳng đứng, vận tốc quay là $\omega = 3000$ vòng/phút. Vòng tròn lớn có bán kính $c = 5.10^{-5}$ (m) bằng khe hở xuyên tâm mô tả phạm vi có thể di chuyển của tâm trục bánh răng. Khi mô phỏng chuyển động của trục, độ lệch tâm $e \geq c$ thì nghĩa là va chạm xảy ra giữa bánh răng và thân bơm. **Hình 5** (phải) cho thấy tâm của trục bánh răng dao động xung quanh vị trí cân bằng. Dao động này mặc dù ngày càng nhỏ nhưng vẫn luôn xảy ra. Trong **Hình 6** (trái) là đồ thị lực tác dụng của màng dầu $F_x$ và $F_y$ theo thời gian. Theo phương pháp phân tích, màng dầu được xem như một hệ lò xo và giảm xóc. Do đó, dao động của trục bơm là dao động tắt dần. Lực $F_x$ hội tụ dần về 0 và lực $F_y$ hội tụ dần về giá trị 1000 (N) cân bằng với tải trọng $W$. Kết quả này hoàn toàn phù hợp với phân tích lý thuyết về chuyển động của trục bơm. **Hình 6** (phải) mô tả phân bố áp suất của lớp màng dầu tại vị trí cân bằng.

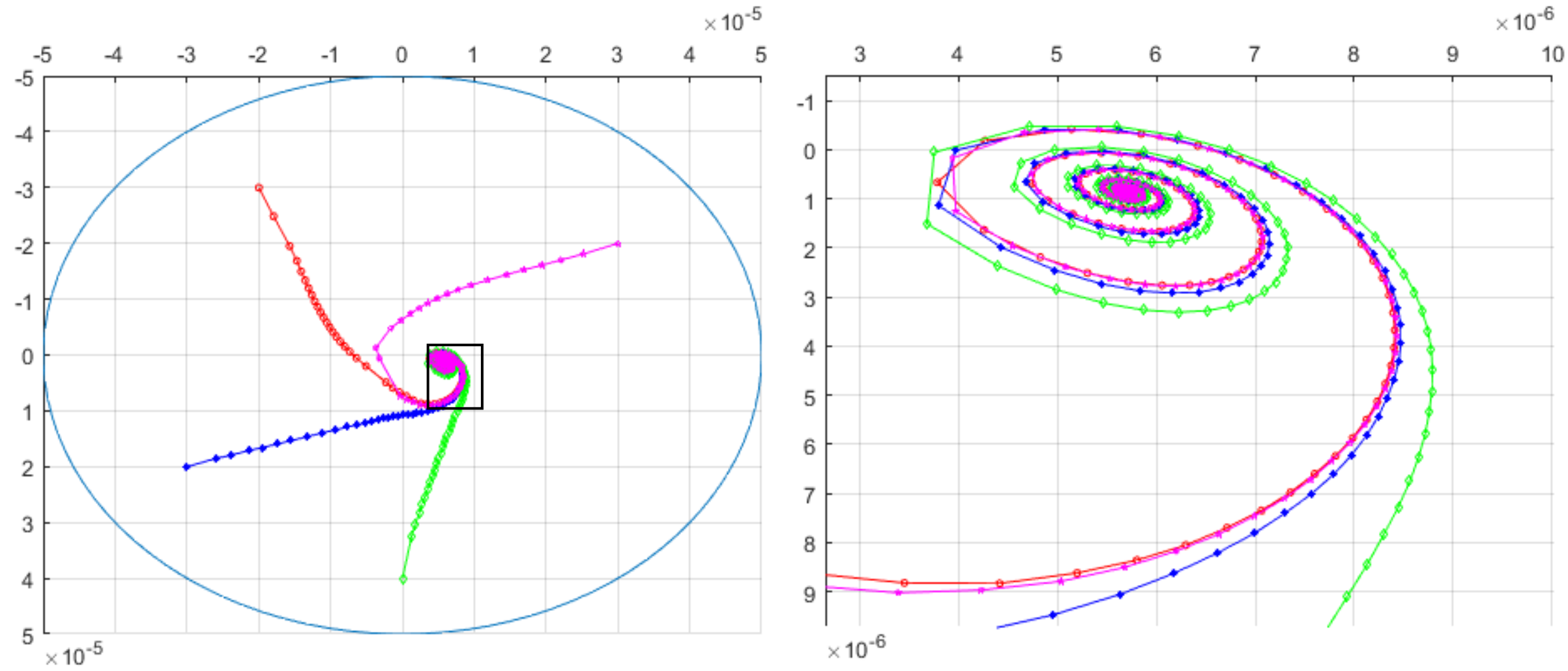

**Hình 7**: Quỹ đạo chuyển động với các vị trí xuất phát khác nhau, $W = 1000$ (N)

**Hình 7** là kết quả mô phỏng quỹ đạo chuyển động của trục bánh răng khi ngoại lực tác dụng là $W = 1000$ (N) và tâm của trục xuất phát tại các vị trí khác nhau. Các quỹ đạo này hội tụ về cùng vị trí cân bằng sau một số lần lặp. Kết quả này là phù hợp với phân tích lý thuyết về chuyển động của trục bơm.

*b) Ảnh hưởng của tải trọng đối với quỹ đạo chuyển động:*

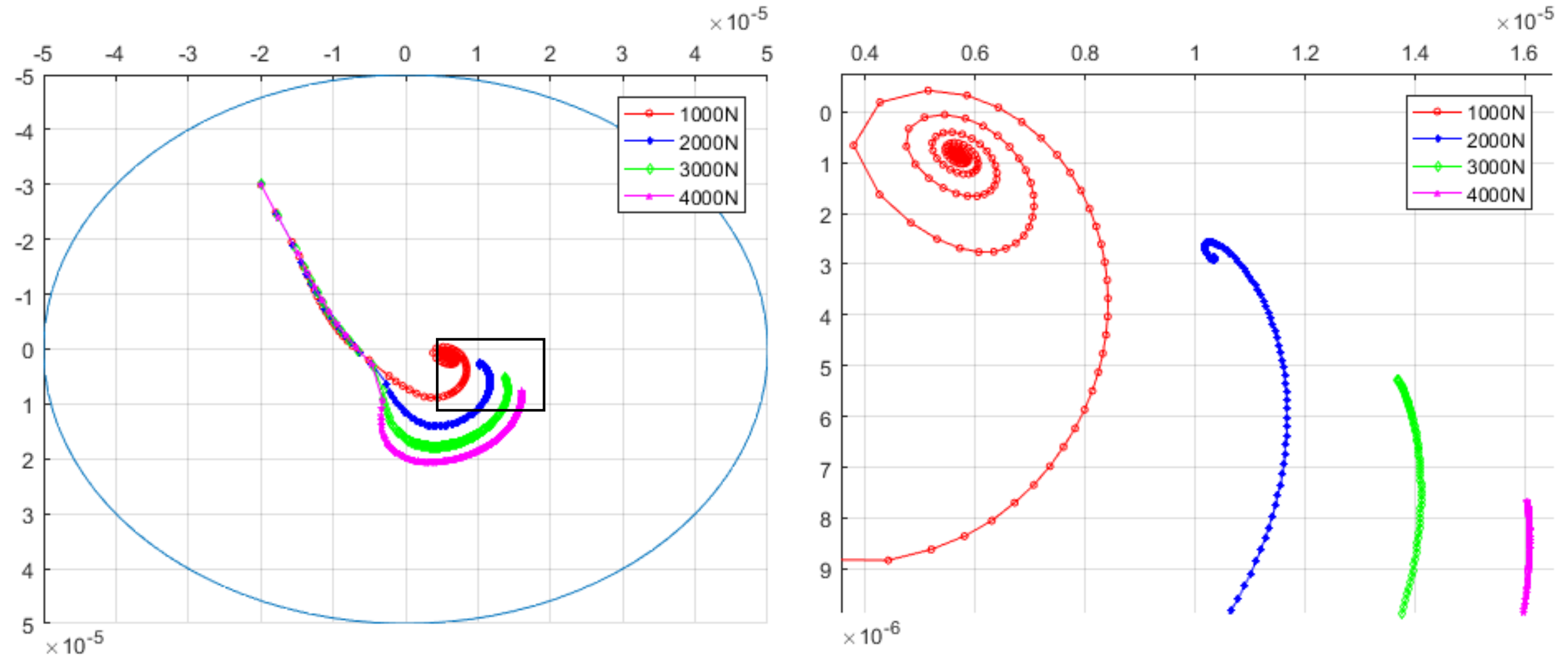

**Hình 8**: Ảnh hưởng của tải trọng đối với chuyển động của trục bơm

Trong **Hình 8**, chúng ta xem xét ảnh hưởng của giá trị tải trọng W đối với chuyển động của trục bơm. Khi ta tăng giá trị của ngoại lực tác dụng từ 1000 N đến 4000 N, kết quả cho thấy trục bơm chuyển động dần về vị trí cân bằng càng thấp (độ lệch tâm càng lớn). Hơn nữa có thể quan sát thấy, khi tải trọng tăng lên thì trục bơm có dao động xung quanh vị trí cân bằng với biên độ càng nhỏ. **Hình 9** biểu diễn mối quan hệ giữa độ lệch tâm và tải trọng tác dụng. Nói chung, tải trọng càng lớn thì độ lệch tâm càng lớn. Tuy nhiên, khi tỉ số độ lệch tâm $\varepsilon = e/c$ tiếp cận 0,8 thì để tăng độ lệch tâm đòi hỏi tải trọng tăng lên khá nhiều.

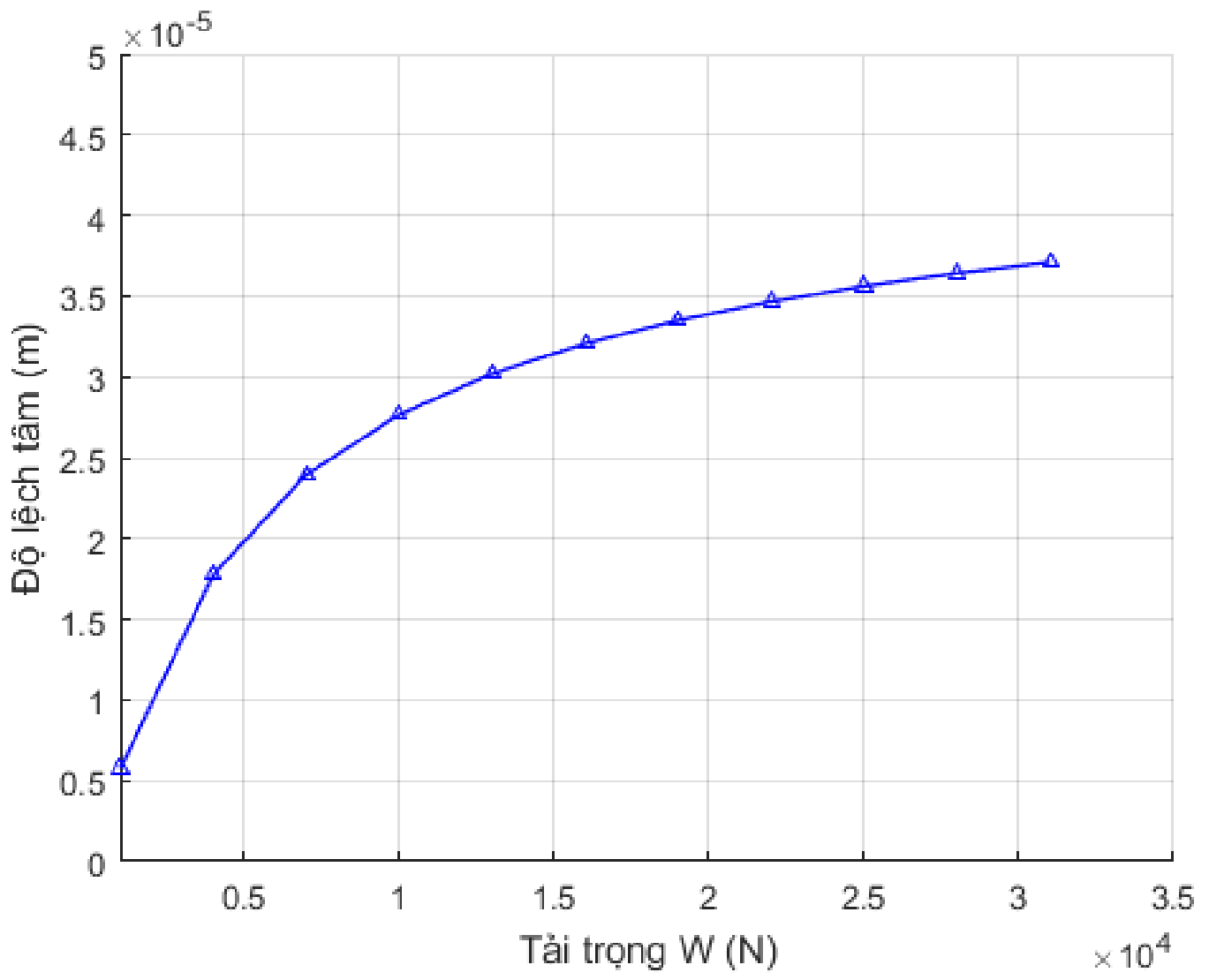

**Hình 9**: Độ lệch tâm theo tải trọng tác dụng

*c) Ảnh hưởng của tốc độ quay đối với quỹ đạo chuyển động:*

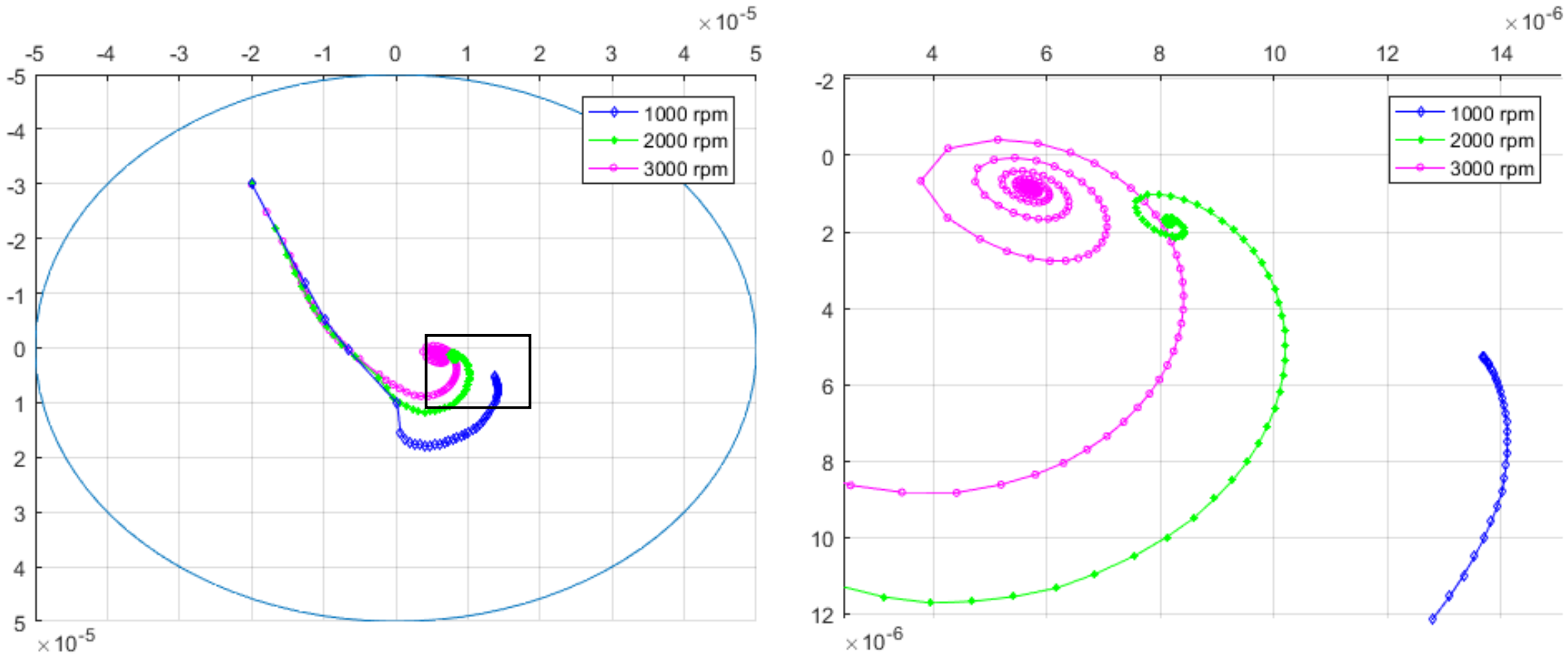

**Hình 10**: Ảnh hưởng của tốc độ quay đối với chuyển động của trục bơm

Khi thay đổi tốc độ quay của trục bánh răng, quỹ đạo chuyển động cũng có sự thay đổi. Cụ thể, vị trí cân bằng lệch tâm càng nhiều khi giảm tốc độ quay từ 3000 vòng/phút xuống 1000 vòng/phút như được mô tả trong **Hình 10**. Đồng thời, trục bánh răng dao động càng lớn quanh vị trí cân bằng khi vận tốc quay càng cao.

*d) Quỹ đạo của trục bánh răng dưới tác động của tải trọng biến đổi:*

Xét trường hợp trục bánh răng chuyển động dưới tác động của tải trọng biến đổi. Cụ thể ta cho thêm vào tải trọng theo phương ngang có dạng hình sin và có biên độ bằng 10% tải trọng theo phương thẳng đứng: $W_x = 10\%.W.\sin \omega t$. Khi đó kết quả mô phỏng về cơ bản giống như trường hợp tải trọng cố định. Tuy nhiên vẫn có điểm khác biệt. Trong trường hợp tải trọng biến đổi, trục bánh răng không chuyển động về vị trí cân bằng mà quay xung quanh vị trí cân bằng như ghi nhận trong **Hình 11** (trái). Điều này xảy ra là do dao động của tải trọng và áp lực màng dầu, xem **Hình 11** (phải).

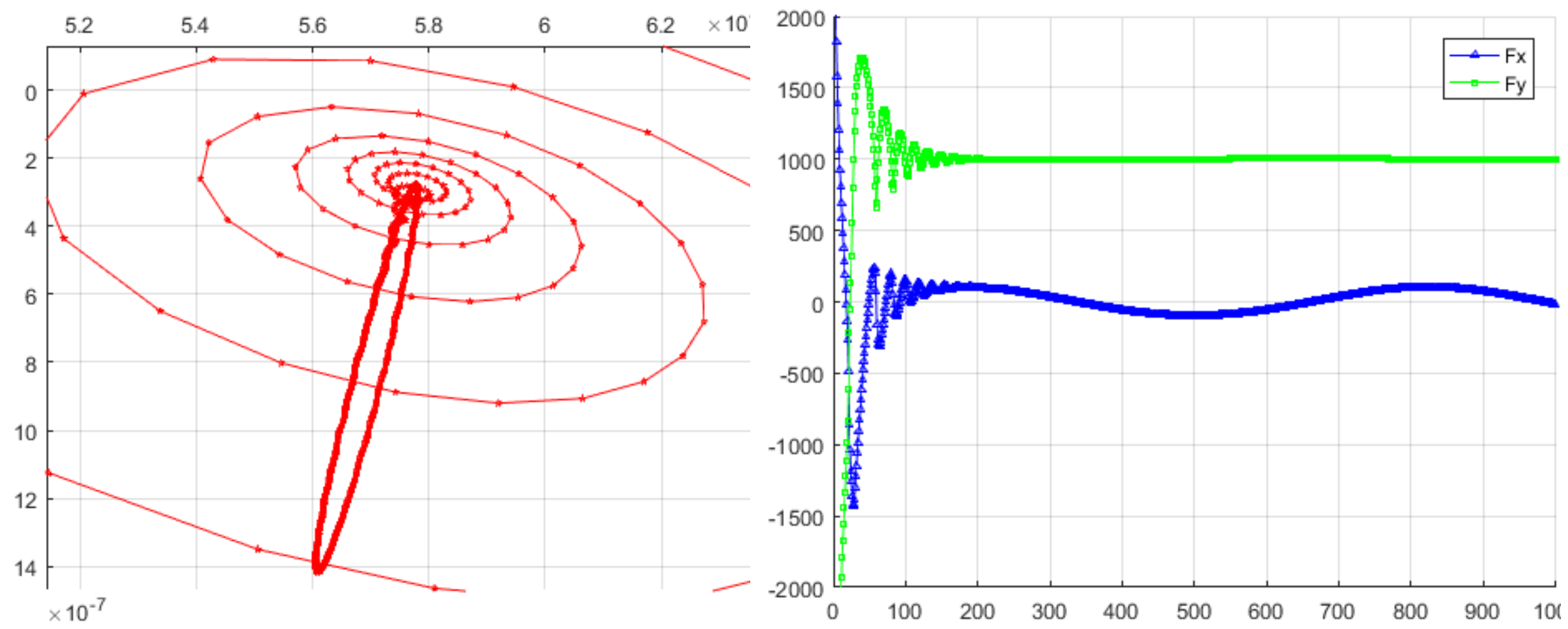

**Hình 11**: Quỹ đạo của trục bơm khi tải trọng thay đổi

# IV. Nhận xét và kết luận

Bài báo này cung cấp một mô hình toán đơn giản, dễ sử dụng để xác định chuyển động của trục trong bơm bánh răng ăn khớp trong. Mô hình là sự kết hợp của các phương trình Reynold mô tả phân bố áp suất màng dầu, các phương trình tính lực, hệ số độ cứng và hệ số giảm xóc của màng dầu, và phương trình chuyển động của trục bánh răng. Thuật toán để tính toán chuyển động được đề xuất và một số kết quả mô phỏng cũng được thực hiện. Các kết quả này có thể làm tiền đề cho việc phân tích trạng thái làm việc của bơm bánh răng ăn khớp trong và xác định các điều kiện tới hạn hay tối ưu để bơm hoạt động.

Trong mỗi bước tính toán mô phỏng đều phải tính toán lại các thông số của mô hình lò xo - giảm chấn. Nói các khác các thông số của mô hình là thay đổi theo các điều kiện làm việc của mỗi lần khảo sát. Kết quả này có được là do việc tuyến tính hóa lực tác động màng dầu, trong thực tế là hàm phi tuyến đối với tọa độ và vận tốc. Việc chia thời gian thành các khoảng nhỏ để đảm bảo quá trình tính toán là ổn định và dẫn đến khối lượng tính toán tăng lên. Đó là sự đánh đổi để thu được một mô hình tính toán đơn giản, dễ cài đặt và có độ chính xác cao. Đối với các máy tính có cấu hình tốt, thời gian tính toán sẽ không trở thành vấn đề lớn.

## Tài liệu tham khảo